\documentclass[twocolumn]{revtex4-1}
\pdfoutput=1
 	 	
\usepackage{subcaption}
\usepackage{enumerate}
\usepackage{booktabs} 
\usepackage{color}
\usepackage{amsfonts}
\usepackage{amsmath}
\usepackage{graphicx}

\usepackage[linesnumbered, ruled]{algorithm2e}

\usepackage{mathtools}
\usepackage{enumerate}
\usepackage{amssymb}
\usepackage{amsthm}

\newcommand*{\trace}{\mathsf{Tr}} 
\DeclareMathOperator{\Span}{span}

\begin{document}

\author{Moritz August}
\affiliation{Department of Informatics, Technical University of Munich, 85748 Garching, Germany (august@in.tum.de)}
\author{Mari Carmen Ba\~{n}uls}
\affiliation{Max Planck Institute for Quantum Optics, 85748 Garching, Germany (mari.banuls@mpq.mpg.de)}
\author{Thomas Huckle}
\affiliation{Department of Informatics, Technical University of Munich, 85748 Garching, Germany (huckle@in.tum.de)}

\title{On the Approximation of Functionals of Very Large Hermitian Matrices represented as Matrix Product Operators}

\begin{abstract}
We present a method to approximate functionals $\trace f(A)$ of very high-dimensional hermitian matrices $A$ represented as Matrix Product Operators (MPOs). Our method is based on a reformulation of a block Lanczos algorithm in tensor network format. We state main properties of the method and show how to adapt the basic Lanczos algorithm to the tensor network formalism to allow for high-dimensional computations. Additionally, we give an analysis of the complexity of our method and provide numerical evidence that it yields good approximations of the entropy of density matrices represented by MPOs while being robust against truncations.
\end{abstract}

\maketitle

\section{Introduction}
Approximating functionals of very large matrices is an important problem in many fields of science, such as network analysis \cite{baglama2014analysis, fenu2013block, estrada2010network, newman2010networks} or quantum mechanics \cite{schollwock2011density, verstraete2008matrix}. In many cases, the respective matrices are hermitian due to either the underlying physical properties of the systems they describe or the way they are constructed from, e.\ g., a graph. Naturally, as the dimensionality of the matrices becomes very high, i.\ e., several tens or hundreds of thousands and above, explicit methods of function evaluation, like exact diagonalization, break down and approximations must be made.

One paradigm for the approximation of high-dimensional matrices that has gained a lot of attention especially in the quantum information, condensed matter and numerical linear algebra communities are tensor network representations \cite{verstraete2008matrix, schollwock2011density, grasedyck2013literature, bachmayr2016tensor, oseledets2011tensor, huckle2013computations}. Among the class of tensor networks, matrix product states (MPS) and matrix product operators (MPO) count among the best established methods. These representations approximate large tensors by contractions of multiple low-rank tensors in a row and have been shown to yield efficient parametrizations of many relevant states of quantum many-body systems~\cite{verstraete2006matrix, perez2006matrix, hastings2007area}.

In this work, we introduce a novel method to approximate functionals of the form $\trace f(A)$ where we assume $f$ to be smooth and defined on the spectrum of $A$ as well as $A$ to be hermitian and given in MPO-format. For our method, we have reformulated a block version of the Lanczos algorithm in MPO/MPS-representation. This particular block Lanczos algorithm will be referred to as global Lanczos algorithm in the following and has already been used to approximate functionals of the given form for explicitly stored matrices \cite{bellalij2015bounding, elbouyahyaoui2009algebraic}. Rewriting it for the tensor network formalism however allows us to consider block vectors of size identical to $A$ which was previously prohibitive. Our method is thus able to approximate $\trace f(A)$ for certain $f(A)$ needing only one carefully selected starting block vector. This means that we get rid of the approximation error induced by having to combine the results obtained for multiple different starting block vectors. At the same time, we find the numerical error induced by the MPS/MPO-representation to be comparably small. Our method can be applied whenever $A$ is efficiently approximated by an MPO. We will in the following however focus on the case where $A$ has directly been defined as MPO.

The rest of this work is structured as follows: after a basic introduction to matrix product states and operators in Section~\ref{MPO}, we will introduce the block Lanczos method we employ in this work in Section~\ref{lanczos}. In Section~\ref{quadrature}, we will then provide an explanation of how one can use the method to approximate functionals of the given type. Following this, we will then state our method in Section~\ref{algo}, show how we have reformulated the global Lanczos method in the tensor network formalism and discuss its properties as well as give an analysis of its complexity. Finally, in Section~\ref{results} we will provide numerical evidence for the fast and robust convergence of our method for the case of the trace-norm and von-Neumann entropy of quantum mechanical density matrices. We conclude with a discussion of the results in Section~\ref{discussion}.
\section{Matrix Product States and Operators}
\label{MPO}
In the area of tensor networks, MPS and MPOs form a well-established class of tensor decompositions that allow for efficient and stable approximation of very high-dimensional vectors and matrices respectively. While they are commonly used in theoretical and numerical quantum many body physics to model, e.g., ground and thermal states \cite{verstraete2004density, zwolak2004mixed, pirvu2010matrix, fannes1992finitely, vidal2003efficient, verstraete2004matrix, schollwock2011density, verstraete2008matrix, verstraete2006matrix}, they also have been independently introduced in the numerical mathematics community under the name of Tensor Trains (TT) \cite{oseledets2011tensor} as a general approximation tool. Since our work is mainly motivated by applications in quantum physics, we will adapt the respective terminology in the following.

A matrix product state is a decomposition of a vector $v \in \mathbb{C}^N$ such that 
\begin{equation*}
v_i = v_{i_1 \dots i_L} = \trace A_1^{i_1}A_2^{i_2} \dotsm A_L^{i_L}
\end{equation*} 
where the index $i$ is split up into $L$ sub-indices of dimensionality $d$, called physical indices, $i_1, \dotsc ,i_L$. $A_1, \dotsc ,A_L \in \mathbb{C}^{d \times D \times D}$ we will refer to as the sites of the MPS and $D$ is called the bond dimension. This concept of splitting up the index $i$ is known under the name of quantized tensor trains (QTT)~\cite{khoromskij2011dlog} in the numerical community. While in principle every site may have its own bond dimensions, as long as they allow for contraction with neighbouring sites, for simplicity and without loss of generality we will assume all sites to have identical bond dimension $D$. The physical dimension $d$ is likewise assumed to be identical for all sites. It is important to note that $N = L^d$, as this relation forms the basis for the ability of MPS/MPOs to represent vectors and matrices of very high dimensionality.

A slightly different representation can be chosen, where 
\begin{equation*}
v_i = v_{i_1 \dots i_L} = A_1^{i_1}A_2^{i_2} \dotsm A_L^{i_L}
\end{equation*}
with $A_1 \in \mathbb{C}^{d \times 1 \times D}$ and $A_L \in \mathbb{C}^{d \times D \times 1}$. In physical terms, the former representation corresponds to systems with closed boundary-conditions (CBC) whereas the latter assumes open boundaries (OBC). It is clear that OBC is just a special case of PBC. For the remainder of this work, we assume open boundaries without loss of generality.
Following this decomposition, a particular element of $v$ is described by a chain of matrix multiplications, explaining the name of the representation.

Now, the whole vector $v$ can be written as 
\begin{align*}
v &= \sum_{i_1,\dots,i_L}^d (A_1^{i_1} A_2^{i_2}\dotsm A_L^{i_L}) \, (e_{i_1}\otimes e_{i_2}\otimes \dotsm \otimes e_{i_L}) \\
&= \sum_{i_1,\dots,i_L}^d \sum_{k_2,\dots,k_{L-1}}^{D} (A_{1,1k_2}^{i_1} A_{2,k_2k_3}^{i_2}\dotsm A_{L,k_{L-1}1}^{i_L}) \, \nonumber \\& 
\qquad \qquad \qquad \qquad \qquad \qquad \cdot (e_{i_1}\otimes e_{i_2}\otimes \dotsm \otimes e_{i_L}) \\
&= \sum_{k_2,\dots,k_{L-1}}^{D} \left(\sum_{i_1}^d A_{1,1k_2}^{i_1}e_{i_1}\right)\otimes \left(\sum_{i_2}^d A_{2,k_2k_3}e_{i_2}^{i_2}\right)\nonumber \\&
\qquad \qquad \qquad \qquad \qquad \otimes \dotsm \otimes \left(\sum_{i_L}^d A_{L,k_{L-1}1}^{i_L}e_{i_L}\right) \\
&= \sum_{k_2,\dots,k_{L-1}}^{D} u_{1,k_2} \otimes u_{2,k_2k_3} \otimes \dotsm \otimes u_{L,k_{L-1}}
\end{align*}
where $e_j$ denotes the $j$th column of the identity matrix. This expression sheds some light on the underlying tensor product structure of MPS and facilitates comparisons with other tensor decomposition schemes.

We now turn our attention to the representation of operators and matrices respectively. Abstractly, one can define an MPO as an operator with an MPS representation in a direct product basis of the operator linear space. More concretely, for the representation of a matrix $M \in \mathbb{C}^{N \times N}$ as an MPO, the approach presented above can easily be adapted to yield
\begin{equation*}
M_{ij} = M_{i_1\dots i_Lj_1\dots j_L} = A_1^{i_1j_1}A_2^{i_2j_2} \dotsm A_L^{i_Lj_L}
\end{equation*}
where $i,j$ have been split up as before and $A_1,\dotsc,A_L \in \mathbb{C}^{d \times d \times D \times D}$. In analogy to the case for a vector, we can write the whole matrix as
\begin{align*}
M &= \sum_{i_1,\dots,i_L}^d \sum_{j_1,\dots,j_L}^d (A_1^{i_1j_1} A_2^{i_2j_2}\dotsm A_L^{i_Lj_L}) \, \nonumber \\& 
\qquad \qquad \cdot (e_{i_1}\otimes e_{i_2}\otimes \dotsm \otimes e_{i_L}) (e_{j_1}^T\otimes e_{j_2}^T\otimes \dotsm \otimes e_{j_L}^T) \\
&= \sum_{i_1,\dots,i_L}^d \sum_{j_1,\dots,j_L}^d \sum_{k_2,\dots,k_{L-1}}^{D} (A_{1,1k_2}^{i_1j_1} A_{2,k_2k_3}^{i_2j_2}\dotsm A_{L,k_{L-1}1}^{i_Lj_L}) \, \nonumber \\& 
\qquad \qquad \qquad \qquad \cdot(e_{i_1}e_{j_1}^T)\otimes (e_{i_2}e_{j_2}^T)\otimes \dotsm \otimes (e_{i_L}e_{j_L}^T) \\
&= \sum_{k_2,\dots,k_{L-1}}^{D} U_{1,k_2} \otimes U_{2,k_2k_3} \otimes \dotsm \otimes U_{L,k_{L-1}}
\end{align*}
where $e_j$ is again the $j$th column of the identity and $e_j^T$ is its transpose. Note that this also holds true for other product bases, like for instance the Pauli basis. Making use of these formulations, it is easy to show that basic operations such as scalar multiplication, addition and inner product as well as the multiplication of an MPS by an MPO or of two MPOs can be performed in the formalism. The addition and non-scalar multiplication however lead to an increase in the bond dimension $D$. For the addition of two MPS/MPOs with bond dimensions $D$ and $D^{\prime}$, the new bond dimension is $D^{\prime \prime}=D+D^{\prime}$ and for the multiplication, $D^{\prime \prime}=D \cdot D^{\prime}$ \cite{schollwock2011density}. This can again easily be verified.

It is obvious from the above explanation that the bond dimension is the decisive factor for the expressive power of the formalism. 
An exact representation of a vector (operator) as an MPS (MPO) is always possible if we allow the bond dimension, $D$, to be big enough, which may mean exponentially large in $L$, up to $d^{N/2}$~\cite{verstraete2004density}. 
When the maximum value of $D$ is limited to some small value (\emph{truncated}) not all vectors or operators can be represented, what may give raise to approximation errors.
We will in the following denote that some vector $v$ or matrix $M$ is approximated with bond dimension $D$ by writing $v[D]$ and $M[D]$ respectively. Nevertheless, it has been found that MPS/MPOs often yield good approximations for $D \in \mathcal{O}(poly(L))$ leading to the total number of parameters $LdD^2 \in \mathcal{O}(poly(L))$ as opposed to $d^L$ or $d^2L$  for the whole vector or matrix, respectively. This constitutes another main reason for their usefulness as an efficient approximation scheme. 

Naturally, many methods have been developed to find optimal and canonical representations for a given $D$ both in the numerical and the quantum physics community. The most important algorithms for optimizing a given MPS/MPO with respect to some error function and bond-dimension thereby rely on local updates of the individual $A_i$, with all other sites being treated as constants, rather than considering all parameters simultaneously. These algorithms, starting with the left- or right-most site, generally sweep back and forth over the chain of sites, updating one site per step, until convergence. As all sites not considered in a given step are treated as fixed, this sweeping scheme allows for reusage of previously computed values in a dynamical programming fashion. As explaining the details and the complexity of these algorithms exceeds the scope of this section, we refer the interested reader to the overview articles \cite{schollwock2011density, verstraete2008matrix, bachmayr2016tensor, grasedyck2013literature}.

\section{The Global Lanczos Algorithm}
\label{lanczos}
The idea of employing variants of the Krylov method to solve various types of problems, for instance solving linear systems \cite{saad1986gmres, meyer1998idea, kressner2011low}, finding eigenvectors \cite{calvetti1994implicitly, arnoldi1951principle, lanczos1950iteration, krylov1931numerical} or approximating the action of an exponential operator onto a vector \cite{garcia2006time}, is already well-established. All the (non-block) methods have in common that, starting with a given matrix $A \in \mathbb{C}^{N \times N}$ and an initial vector $u \in \mathbb{C}^N$, they construct an orthonormal basis $U_{K} = \left[ u_1, u_2,\dotsc,u_K \right] \in \mathbb{C}^{N \times K}$ of the Krylov space $\mathcal{K}\mathcal{R}_{K}=\{u, Au, A^2u,\dotsc,A^{K-1}u\}$, where $K$ is the dimension of the space. We will now first describe the original Lanczos method and then generalize our findings to the global algorithm.
\begin{algorithm}[t]
\setlength{\leftskip}{10pt}
\setlength{\skiprule}{10pt}
\caption{Original Lanczos Algorithm}\label{vector_lanczos}
    \SetKwInOut{Input}{Input}
    \SetKwInOut{Output}{Output}

    \Input{Matrix $A \in \mathbb{C}^{N \times N}$, Starting Vector $u \in \mathbb{C}^N$, Number of Dimensions $K$}        
    $u_0 \leftarrow 0$ \;
    $v_0 \leftarrow u $ \;
    \For{$i \leftarrow 1 ; i \leq K $}{
	$\beta_i \leftarrow \|v_{i-1}\|$ \;
    \If{$\beta_i = 0$}{\quad break \;}    
    $u_i \leftarrow  v_{i-1} / \beta_i$ \;
    $v_i \leftarrow Au_{i} - \beta_i u_{i-1}$ \;
    $\alpha_i \leftarrow u_i^*v_i$ \;
    $v_i \leftarrow v_i - \alpha_i u_{i}$ \;    
    }
    
  	\Output{Orthonormal Basis $U_{K} \in \mathbb{C}^{N \times K}$, Tridiagonal Matrix $T_K \in \mathbb{R}^{K \times K}$}
\end{algorithm}
The original Lanczos algorithm is a variant of the general Krylov method that assumes $A$ to be hermitian and in our case makes use of the Gram-Schmidt (GS) orthogonalization method to construct $U_K$. During the process of building $U_K$, the Lanczos method produces a matrix $T_K \in \mathbb{R}^{K \times K}$ that is given by
\begin{equation*}
T_K = \left[ \begin{matrix}
\alpha_1 & \beta_2 & & \mathbf{0}\\
\beta_2 & \alpha_2 & \ddots & \\
& \ddots & \ddots & \beta_K \\
\mathbf{0} & & \beta_K & \alpha_K \end{matrix} \right],
\end{equation*}
where $\alpha_i$ and $\beta_i$ are defined as in Algorithm \ref{vector_lanczos}. Reminding ourselves of the orthonormality of $U_K$ and the hermiticity of $A$, it is easy to verify that
\begin{equation*}
U_K^* A U_K = T_K,
\end{equation*}
where $U^*$ represents the Hermitian conjugate, and that hence $A$ is similar to $T_N$, $T_N$ being the tridiagonal matrix of full dimension $N$. This provides an interesting perspective on the Lanczos method as a projection onto a lower-dimensional space and a first glance at how one may approximate the spectrum of $A$ via the eigenvalues of $T_K$, called the Ritz values of $A$ \cite{saad1992numerical, trefethen1997numerical}. Note that $T_K$ has to be real because the $\beta_i$ are vector-norms and the $\alpha_i$ have to be real due to the hermiticity of $A$ and its similarity to $T_N$. We now define an extended version $\widehat{T}_K$ of $T_K$ by 
\begin{equation*}
\widehat{T}_K = \left[ \begin{matrix}
\alpha_1 & \beta_2 & & \mathbf{0}\\
\beta_2 & \alpha_2 & \ddots & \\
& \ddots & \ddots & \beta_K \\
\mathbf{0} & & \beta_K & \alpha_K \\
0&\dots&0&\beta_{K+1} \end{matrix} \right]
\end{equation*}
and use it to state the partial Lanczos decomposition
\begin{equation*}
A U_K = U_{K+1} \widehat{T}_K = U_K T_K + \beta_{K+1}u_{K+1}e_K^T
\end{equation*} that links $U_K$ to $U_{K+1}$ for $K<N$. The decomposition makes it clear that the improvement from one step of the algorithm to the next is proportional to $\beta_{K+1}$. In the case where $\beta_{K+1} = 0$, the decomposition is exact and $U_K$ spans an eigenspace of $A$. For the decomposition, we have implicitly used a three-term recurrence relation that can be read off directly from Algorithm~\ref{vector_lanczos}. By construction it holds that
\begin{equation*}
Au_i = \beta_{i+1} u_{i+1} + \beta_i u_{i-1} + \alpha_i u_i
\end{equation*}
which can be rewritten to yield
\begin{equation}
\label{lanczos_rec}
\beta_{i+1} u_{i+1} = (A - \alpha_i I_N) u_i - \beta_i u_{i-1}
\end{equation}
where $I_N$ is the $N \times N$ identity matrix.
This relation plays an important role in our algorithm and we will get back to it when we explain how we approximate functionals in the next section.

However, before explaining the rigorous approach, we will provide two more intuitive arguments to show why Krylov methods can be used for our purpose. From the above explanations, it is evident that every vector $v \in \Span{(U_K)}$ is given by 
\begin{equation*}
v = \sum_{i=0}^{K-1} c_i(A^iu) = \left(\sum_{i=0}^{K-1} c_iA^i \right) u = p_v(A) u
\end{equation*}
where $p_v(A)$ is some polynomial of $A$ with degree $K-1$. Note that there exist instances of $v$ for which $p_v(A)$ yields a power-series approximation around 0 of degree $K-1$ of some desired function $f(A)$. This is a key insight if one wants to understand why Krylov methods are very useful for problems like the ones mentioned above. It is obvious from this formulation that $K$ has a strong impact on the error of the approximation of $f(A)$ which naturally leads to asking whether bounds can be derived for $K$ with respect to certain functions. It can be seen that the $K$ needed for exactness can be determined from the degree of the minimal polynomial of $A$ \cite{meyer1998idea}. We have seen above that $U_K^* A U_K = T_K$. Therefore it holds that 
\begin{align*}
U_K^* f(A) U_K & \approx U_K^* \left( \sum_{i=0}^{K-1} c_i A^i \right) U_K & \\\ & \approx \sum_{i=0}^{K-1} c_i (U_K^* A U_K)^i \\ & =  \sum_{i=0}^{K-1} c_i  T_K^i \\ & \approx f(T_K),
\end{align*} where we again assume some power-series approximation of $f(A)$ of degree $K-1$ and make use of the fact that $\lim_{K \rightarrow N} U_K^*U_K = I$.

After having reviewed the original algorithm, we now proceed to state the global block version used in this work. In order to do so, we will from now on assume to be given an \emph{initial matrix} $U \in \mathbb{C}^{N \times M}$. Starting from $U$, the algorithm will then build up a basis of matrices $\mathbf{U}_K = [U_1,\cdots,U_K]$ with $\mathbf{U}_K \in \mathbb{C}^{N \times KM}$. While there exist several block-versions of the Lanczos algorithm \cite{elbouyahyaoui2009algebraic, grimes1994shifted, montgomery1995block, golub1981block, cullum1974block}, we will only consider the one presented in \cite{bellalij2015bounding} as it does not require the columns of $U_i$ to be orthogonal, which would be prohibitive for very large matrices.
\begin{algorithm}[t]
\setlength{\leftskip}{10pt}
\setlength{\skiprule}{10pt}
\caption{Global Lanczos Algorithm}\label{global_lanczos}
    \SetKwInOut{Input}{Input}
    \SetKwInOut{Output}{Output}

    \Input{Matrix $A \in \mathbb{C}^{N \times N}$, Starting Matrix $U \in \mathbb{C}^{N \times M}$, Number of Dimensions $K$}    
    $U_0 \leftarrow 0$ \;
    $V_0 \leftarrow U$ \;
    \For{$i \leftarrow 1 ; i \leq K $}{
    $\beta_i \leftarrow \|V_{i-1}\|_F$ \;
    \If{$\beta_i = 0$}{break \;}
    $U_i \leftarrow  V_{i-1} / \beta_i$ \;
    $V_i \leftarrow AU_{i} - \beta_i U_{i-1}$ \;
    $\alpha_i \leftarrow \langle U_i, V_i \rangle$ \;
    $V_i \leftarrow V_i - \alpha_i U_{i}$ \;
    }
    
  	\Output{Orthonormal Basis $\mathbf{U}_{K} \in \mathbb{C}^{N \times KM}$, Tridiagonal Matrix $T_K \in \mathbb{R}^{KM \times KM}$}
\end{algorithm} Now, we first need to state the inner product with respect to which the individual $U_i$ must be orthonormal and define it to be
\begin{equation*}
\langle U_i, U_j \rangle = \trace U_i^* U_j
\end{equation*}
where $U_i,U_j \in \mathbb{C}^{N \times M}$. This induces the well known Frobenius norm
\begin{equation*}
\| U_i \|_F = \langle U_i, U_i \rangle^{1/2}
\end{equation*}
and hence this definition of the inner product is a straight-forward generalization of the one used in Algorithm \ref{vector_lanczos}. Naturally, one may also choose different inner products \cite{elbouyahyaoui2009algebraic}. For this work, we do however choose the Frobenius norm as it can be efficiently computed for MPOs. Equipped with this definition, we can see that Algorithm \ref{global_lanczos} is in fact a direct generalization to the matrix-case. As such we find that after $K$ steps, the method has produced the reduction of $A$ $T_K$ and yields the partial global Lanczos decomposition
\begin{equation*}
A \mathbf{U}_K = \mathbf{U}_K \tilde{T}_K + \beta_{K+1}U_{K+1}E_K^T
\end{equation*}
where we define $\tilde{T}_K = T_K \otimes I_M \in \mathbb{R}^{KM \times KM}$ and $E_K^T = [\mathbf{0},\cdots,\mathbf{0},I_M] \in \mathbb{R}^{M \times KM}$. Furthermore, it again holds that
\begin{equation}
\label{lanczos_rec_glob}
\beta_{i+1} U_{i+1} = (A - \alpha_i I_N) U_i - \beta_i U_{i-1}
\end{equation}
and
\begin{equation*}
\mathbf{U}_K^* A \mathbf{U}_K = T_K
\end{equation*}
if we apply the inner product defined previously. From now on, we will implicitly make use of this inner product whenever appropriate. Then, all other observations made for the original Lanczos method above carry over to the global Lanczos case.
\section{Introducing Gauss-type Quadrature}
\label{quadrature}
In the previous section, we have sketched that for the Lanczos method $U_K^* f(A) U_K \approx f(T_K)$ and hence \begin{equation*}
u_1^* f(A) u_1 = e_1^T U_K^* f(A) U_K e_1 \approx e_1^T f(T_D) e_1,
\end{equation*}
which generalizes to the global Lanczos method. However, there exists a well studied and more rigorous approach to function approximation via the Lanczos methods \cite{bai1996some, golub2009matrices, golub1994matrices}. This approach revolves around the connection between $T_K$ and Gauss-type quadrature rules. 

To establish the link between both methods, we start by observing that \begin{align*}
u_1^* f(A) u_1 & = u_1^* V f(\Lambda) V^* u_1 \\ & = \sum_{i=1}^N f(\lambda_i) \mu_i^2 \\ & = \int_a^b f(\lambda) \mathrm{d}\mu(\lambda)
\end{align*} 
with $V \Lambda V^*$ being the spectral decomposition of $A$, $\mu_i = e_i^T V^* u_1$ and \begin{equation*}
\mu(\lambda) = 
\begin{cases}
0 & \quad \text{if } \lambda < \lambda_1 = a \\
\sum_{i=1}^j \mu_i^2 & \quad \text{if } \lambda_j \leq \lambda < \lambda_{j+1}  \\
\sum_{i=1}^N \mu_i^2 & \quad \text{if } b = \lambda_N \leq \lambda
\end{cases}
\end{equation*} being a piecewise constant and nondecreasing distribution function. Here we assume the eigenvalues of $A$ to be ordered ascendingly. We can use this result to obtain \begin{equation} \label{integral_stielt}
\begin{split}
\mathcal{I}f &\coloneqq \trace (U_1^* f(A) U_1) \\
			 &= \sum_{i=1}^N e_i^* U_1^* V f(\Lambda) V^* U_1 e_i \\
			 &= \sum_{i=1}^N \int_a^b f(\lambda) \mathrm{d}\mu_i(\lambda) \\
			 &= \int_a^b f(\lambda) \mathrm{d} \sum_{i=1}^N \mu_i(\lambda) \\
			 &= \int_a^b f(\lambda) \mathrm{d}\mu(\lambda)
\end{split}
\end{equation} where we define $\mu_i(\lambda)$ analogously to the case above and $\mu(\lambda) \coloneqq \sum_{i=1}^N \mu_i(\lambda)$. This Riemann-Stieltjes integral can now be tackled via Gauss-type quadrature.

The most general formulation of a Gauss-type quadrature rule is \begin{equation*}
\mathcal{G}f \coloneqq \sum_{i=1}^K \omega_i f(\theta_i) + \sum_{j=1}^M \nu_j f(\tau_j)
\end{equation*} where $\theta_i$ and $\tau_j$ are called the nodes and $\omega_i$ and $\nu_j$ the weights of the quadrature. The remainder of such an approximation is given by \begin{equation*}
\begin{split}
\mathcal{R}f &\coloneqq \int_a^b f(\lambda) \mathrm{d}\mu(\lambda) - \mathcal{G}f \\
             &= \frac{f^{2K+M}(\eta)}{(2K+M)!} \int_a^b \prod_{i=1}^M (\lambda - \tau_i) \left( \prod_{j=1}^K (\lambda - \theta_j) \right)^2 \mathrm{d} \mu (\lambda)
\end{split}
\end{equation*} where $\lambda_1 < \eta < \lambda_N$. Setting $M=0$ yields the Gauss quadrature, while using prescribed nodes results, e. g., in the Gauss-Radau and Gauss-Lobatto quadratures \cite{davis2007methods, bellalij2015bounding}. For the Gauss quadrature, the sign of the remainder thus is given by the sign of the $2K$-th derivative of $f$. Therefore, it is easy to determine whether, for a given $f$, the Gauss quadrature provides a lower or upper bound of the correct value. By application of the Weierstrass-theorem, it can also be shown that \begin{equation*}
\lim_{K \rightarrow \infty} \mathcal{G}_Kf = \mathcal{I}f
\end{equation*} where $\mathcal{G}_Kf$ is the Gauss quadrature with $K$ nodes. Furthermore, $\mathcal{G}_Kf$ is exact for all polynomials of degree at most $2K-1$ \cite{davis2007methods}. It is well known that in order to determine the $\omega_i$ and $\theta_i$ that satisfy this property, one can construct a sequence of polynomials $\{p_0,\cdots,p_K\}$ that are orthonormal in the sense that \begin{equation*}
\int_a^b p_i(\lambda) p_j(\lambda) \mathrm{d}\mu(\lambda) = \delta_{ij}
\end{equation*}
and that satisfy a recurrence relation given by \begin{equation}
\label{gauss_reccur}
\beta_j p_j(\lambda) = (\lambda - \alpha_{j-1}) p_{j-1}(\lambda) - \beta_{j-1} p_{j-2} (\lambda),
\end{equation}
where $p_{-1}(\lambda) \equiv 0$ and $p_0(\lambda) \equiv 1$. Then, the roots of $p_K$ ca be shown to be the optimal $\theta_i$~\cite{davis2007methods, golub2009matrices}. Now, the recurrence relation yields a recurrence matrix $T_K$ defined by \begin{equation*}
T_K = \left[ \begin{matrix}
\alpha_1 & \beta_2 & & \mathbf{0}\\
\beta_2 & \alpha_2 & \ddots & \\
& \ddots & \ddots & \beta_K \\
\mathbf{0} & & \beta_K & \alpha_K \end{matrix} \right]
\end{equation*} whose eigenvalues are the zeroes of $p_K(\lambda)$ and consequentially the $\theta_i$ of $\mathcal{G}f$ \cite{golub2009matrices}. The $\omega_i$ are given by the squared first elements of the normalized eigenvectors of $T_K$ and so, \begin{equation*}
\mathcal{G}f = e_1^T f(T_K) e_1 = e_1^T V f(\Lambda) V^* e_1,
\end{equation*}
where $V \Lambda V^*$ is the spectral decomposition of $T_K$.

Now, the $U_i$ from the global Lanczos method can be expressed by \begin{equation*}
U_i = p_{i-1}(A)U
\end{equation*} with $p_{i-1}$ being some polynomial of degree $i-1$. Then, it is clear that \begin{equation*}
\langle p_{i-1}(A)U, p_{j-1}(A)U \rangle = \langle U_i, U_j \rangle = \delta_{ij}
\end{equation*} and taking into account the above derivations 
\begin{align*}
\langle p_{i-1}(A)U, p_{j-1}(A)U \rangle & = \trace (U^* p_{i-1}(A)^* p_{j-1}(A) U) \\ & = \int_a^b p_{i-1}(\lambda) p_{j-1}(\lambda) \mathrm{d} \mu (\lambda).
\end{align*} Hence, the global Lanczos method produces orthonormal polynomials \cite{lanczos1950iteration} that in addition satisfy the recurrence relation stated in Equation~\ref{gauss_reccur} by construction, as we have shown in Section~\ref{lanczos}. The $T_K$ obtained by the global Lanczos algorithm is thus the recurrence matrix needed to perform the Gauss quadrature. If we choose $U \in \mathbb{C}^{N \times N}$ to be orthonormal, it follows that \begin{equation*}
\trace f(A) = \trace (U^* f(A) U) = \int_a^b f(\lambda) \mathrm{d}\mu(\lambda) \approx e_1^T f(T_K) e_1.
\end{equation*}
\section{Assembling the Parts}
\label{algo}
\begin{algorithm}[t!]
\setlength{\leftskip}{10pt}
\setlength{\skiprule}{10pt}
\caption{Approximation Algorithm}\label{approx_lanczos}
    \SetKwInOut{Input}{Input}
    \SetKwInOut{Output}{Output}
    \SetKwFunction{multOpt}{multiplyAndOptimize}
    \SetKwFunction{addOpt}{sumAndOptimize}
    \SetKwFunction{contract}{contract}
	\SetKwFunction{scalMult}{multiplyScalar}
	\SetKwFunction{spec}{spectralDecomposition}
	\SetKwFunction{checks}{checkStop}

    \Input{MPO $A[D_{A}] \in \mathbb{C}^{N \times N}$, Starting orthogonal MPO $U[D_{init}] \in \mathbb{C}^{N \times N}$, Number of Dimensions $K$, Maximal Bond-Dimension $D_{max}$, Stopping Criteria $\mathcal{S}$}    
    $U_0 \leftarrow 0$ \;
    $V_0 \leftarrow U$ \;
    $D \leftarrow D_{init}$ \;
    \For{$i \leftarrow 1 ; i \leq K $}{
	$\beta_i \leftarrow \sqrt{\contract (V_{i-1}, V_{i-1})}$ \;
    \If{$\beta_i = 0$}{break \;}
    $U_i \leftarrow  \scalMult(1 / \beta_i,V_{i-1})$ \;    
    $D \leftarrow \min(D_{max}, D \cdot D_{A})$ \;
    $V_i \leftarrow \multOpt(A, U_{i}, D)$ \;
	$D \leftarrow \min(D_{max}, D + D_{U_{i-1}})$ \;    
    $V_i \leftarrow \addOpt(V_i, - \beta_i U_{i-1}, D)$ \;
    $\alpha_i \leftarrow \contract(U_i, V_i)$ \;
	$D \leftarrow \min(D_{max}, D + D_{U_{i}})$ \;    
    $V_i \leftarrow \addOpt(V_i, -\alpha_i U_{i}, D)$ \;
    
    $V \Lambda V^* \leftarrow \spec(T_i)$ \;
    $\mathcal{G}f \leftarrow \beta_1^2 e_1^T V f(\Lambda) V^* e_1$ \;
    \If{$\checks(\mathcal{G}f, \Lambda, \mathcal{S})$}{break \;}
    }
    
  	\Output{Approximation $\mathcal{G}f$ of $\trace f(A)$}
\end{algorithm}
Now that we have reviewed the relevant theoretical aspects, we will proceed by showing how we put together the pieces to obtain our algorithm. The whole algorithm is shown in Algorithm \ref{approx_lanczos}. 

Since the global Lanczos method is based on matrix-matrix multiplications, additions of matrices and multiplications of matrices by scalars, these operations have to be formulated for the MPO-case. As the bond dimension of the basis-MPOs grows with the number of multiplications and additions, we need to keep track of the bond dimensions and perform projections onto lower bond dimensions whenever necessary. Thus, the input parameters of our method are the MPO $A[D_A] \in \mathcal{C}^{N \times N}$, an orthogonal starting MPO $U[D_{init}] \in \mathcal{C}^{N \times N}$, the maximal Krylov-dimension $K$, a set of stopping criteria $\mathcal{S}$ and finally the maximal bond-dimension $D_{max}$ of the $U_i$. Note that we assume $U[D_{init}]$ to be orthogonal and of the same dimension as $A[D_A]$. This allows us to replace the approximation that had to be made previously by combining the estimations for several starting matrices by the exact computation. This is only possible due to the fact that we translate the Lanczos method to the MPO-formalism. Besides the case of very large matrices that can be explicitly stored but are too large for multiplications with equally-sized matrices, this is especially important for the case where the respective matrices are only given in MPO-format and $N$ is of the order of several millions, as in the case of quantum many body systems. While we introduce some approximation error by using MPOs, we will show in Section~\ref{results} that these errors can be comparably small already for low bond dimensions in cases of practical interest. In the following, we will omit the declaration of the bond-dimension of an MPO whenever it increases clarity.

While in principle every orthogonal MPO can serve as a starting point, in this work we choose $U[D_{init}]$ to be the identity matrix because it has a minimal MPO-formulation of $A^{jk}_i = \delta_{jk}$.
This allows us to start from the minimal bond-dimension $D_{init}=1$ and thus maximizes the amount of relevant information we can store for a given $D_{max}$. In certain cases it might however be possible to choose a better starting MPO, e. g., when $A$ is very close to being diagonal. For the implementation of the inner product and norm used in the global Lanczos algorithm, we observe that \begin{equation*}
\langle U_i, U_j \rangle = \sum_{k=1}^N \sum_{l=1}^N U_{i,kl}U_{j,lk} = U_{i, vec}^* U_{j,vec}
\end{equation*}
where $U_{i,vec}$ and $U_{j,vec}$ are the vectorized versions of $U_i$ and $U_j$ respectively. This allows us to make use of an efficient, exact way of computing the inner product of MPS \cite{schollwock2011density} by rewriting the $A_i^{j_ik_i}$ as $A_i^{j_i^{\prime}}$ with $\dim j_i^{\prime}=\dim j_i \cdot \dim k_i$ and hence vectorizing the MPO. This functionality is implemented in \texttt{contract}(). The implementation of the scalar multiplication \texttt{multiplyScalar}() is straightforward as it corresponds to multiplying an arbitrary $A_i$, we choose $A_1$ for simplicity, of the respective MPO by the scalar at hand.

A bit more care has to be taken when implementing the functions \texttt{multiplyAndOptimize}() and \texttt{sumAndOptimize}(), the multiplication and summation of MPOs respectively. One possiblity would be to first perform the respective operation exactly, i.e., use the bond-dimension required for the exact result, and to project the resulting MPO onto the current $D$ via the singular value decomposition (SVD) of its $A_i$ in a second step. It has however been found that performing the projection simultaneously to the multiplication or summation at the level of the individual $A_i$ yields superiour results, see \cite{verstraete2008matrix, schollwock2011density, waldherr2014numerical}. In case of the multiplication, we implement this strategy by solving the optimization problem \begin{equation*}
\min_{\tilde{A_i}} \|AU_j[D_{old}] - \tilde{U_j}[D_{new}]\|^2_F
\end{equation*} where $D_{old}$ is the bond dimension used previous to the multiplication and $D_{new}$ is the bond dimension used for the optimization. $\tilde{U_j}[D_{new}]$ denotes the result of the multiplication of $A$ on $U_j$ and the $\tilde{A_i}$ are its tensors. The implementation hence performs the multiplication of the MPOs at tensor level and directly optimizes the resulting tensors for the chosen bond dimension by employing the sweeping scheme sketched in Section \ref{MPO}. In order to apply this algorithm to the case of MPO-MPO multiplication, we rewrite $AU_j$ as \begin{equation*}
(I \otimes A)U_{j,vec} = \left[ \begin{matrix}
A & 0 & & \mathbf{0}\\
0 & A & \ddots & \\
& \ddots & \ddots & 0 \\
\mathbf{0} & & 0 & A \end{matrix} \right]
\left[ \begin{matrix}
U_{j,1}\\
U_{j,2}\\
\vdots\\
U_{j,N} \end{matrix} \right]
\end{equation*} with $U_{j,k}$ being the $k$th column vector of $U_j$. Due to technical reasons, we do in fact use $U_j \otimes I$ and $A_{vec}$. For the summation, we apply the same strategy and solve \begin{equation*}
\min_{\tilde{A_i}} \|\left(U_j[D_{old}] + \sum_k \gamma_k U_k[D_{old}^k]\right) - \tilde{U_j}[D_{new}]\|^2_F
\end{equation*} where $D_{old}$ is the bond dimension used before the addition, $D_{old}^k$ are some other previously used bond dimensions, $D_{new}$ is the bond-dimension to be used for the optimization and $\gamma_k \in \mathbb{C}$ are some scalars. $\tilde{U_j}[D_{new}]$, similar to before, represents the outcome of the summation and the $\tilde{A_i}$ are its tensors.

As it can be seen in Algorithm \ref{approx_lanczos}, we allow for exact multiplication and summation as long as the resulting bond-dimension does not grow beyond $D_{max}$. This however happens quickly, since $D$ grows exponentially with the number of iterations, and so, most of the $U_i$ will be represented based on $D_{max}$. This underlines the importance of $D_{max}$ for the accuracy of the approximation.

After the algorithm has completed one iteration of the global Krylov method, the spectral decomposition of $T_i$ is performed and the current approximation is computed. Based on the approximation and the eigenvalues of $T_i$ the algorithm then determines if it should be stopped in \texttt{checkStop}. Here we have to account for several factors. 

Firstly, we know that $\mathcal{G}f$ converges to the correct value in absence of approximation errors. So, the algorithm can terminate when the distance between the previous and current $\mathcal{G}f$ becomes smaller than some $\epsilon$. 

Secondly, it is clear that the projection of the generated MPOs down to $D_{max}$ introduces an approximation error. While it is possible to obtain the error of the optimization problems described above, it is not clear how the accumulated error influences $\mathcal{G}f$ precisely. However, a possible way of detecting when the approximation error has become too large is to check for the violation of some theoretical constraints. For instance, in case of a positive $A$, we know that the Ritz-values $A$ must be positive as well. If $T_i$ starts to have negative eigenvalues, we thus know that the total approximation error has reached a level that leads to unreasonable results. The same reasoning could be applied for other intervals in which one knows the spectrum of $A$ to be in. 

We have also seen in the previous section that, depending on the sign of the derivative of $f$ in $(\lambda_1, \lambda_N)$, $\mathcal{G}f$ yields an upper/lower bound and that it converges to the true value. Based on this it is possible to show that $\mathcal{G}_M < \mathcal{G}_{M+1}$ for the lower-bound case and $\mathcal{G}_M > \mathcal{G}_{M+1}$ for the case of an upper bound \cite{golub2009matrices}. This provides another stopping-criterion.


As the accumulation of truncation errors can lead to unreasonable results even before the violation of the above property, we propose to keep a moving average over the last $k$ approximations and employ the 3$\sigma$-rule to detect improbable results. This heuristic is justified by the guaranteed convergence in the absence of numerical errors.

After having explained the algorithm, a few remarks are in order:\begin{itemize}

\item In this version of the algorithm, we only consider the Gauss quadrature. This is mainly due to the fact that obtaining good lower or upper bounds on the spectrum of $A$ is in general not possible because of the size of its dimensions. Analogously to \cite{bellalij2015bounding}, our algorithm can nevertheless be adapted to perform Gauss-Radau or Gauss-Lobatto approximations.

\item To improve numerical stability and prevent 'ghost' eigenvalues from occuring, it could be beneficial to perform reorthogonalization. Due to the MPO-representation, this would however be very costly and not necessarily result in a large improvement. Thus, we do not consider this extension. It can however be easily added to the algorithm.

\item In the presented algorithm, we stick to the canonical way of orthogonalizing the new matrix against the old matrices individually. In the case of exactly stored matrices/vectors, this scheme increases numerical stability. Since we now employ approximations of the exact matrices, it might however be worth considering to compute $\alpha_i$ first and then optimize the sum containing both $U_i$ and $U_{i-1}$. The advantage of being able to optimize the whole expression at once might outweigh the disadvantage of orthogonalizing against both matrices simultaneously. On the other hand side, computing $\alpha_i$ first might lead to different and possibly worse results.

\item As we have stated above, it is possible to obtain approximation errors from both \texttt{multiplyAndOptimize} and \texttt{sumAndOptimize}. But these errors naturally only refer to the current optimization and do not allow for strict bounds on the overall error. One could of course try to increase the bond dimension for each individual optimization until its error converges to make sure the partial result is close to exact. The problem here is that due to the possibly exponential growth of the bond dimension needed for exactness, $D_{max}$ is typically reached within very few iterations. From this point on, it is not possible any more to increase $D$ and so, the information about the error provides little useful information. This is why we have resorted to the approach of checking for the violation of theoretically guaranteed behaviour.

\item From the above explanations it is clear that $K$ and $D_{max}$ are the parameters that control the accuracy of the approximation. For the algorithm to be of use for very high-dimensional matrices, we must impose the restriction that $K, D_{max} \in \mathcal{O}(poly(L))$. This property is particularly relevant for quantum mechanical simulations where $N$ grows exponentially with the number of particles.
\end{itemize}

\begin{table}[t]
\centering
\caption{A listing of the complexity of the subfunctions of Algorithm \ref{approx_lanczos}. $L = \log (N)$ is the number of tensors of the MPOs, $d$ is the physical dimension. For simplicity, all $U_i$ are assumed to have bond-dimension $D_{max}$ and $T_i$ is assumed to be in $\mathbb{R}^{K \times K}$.}
\renewcommand{\arraystretch}{1.5}
\begin{tabular}{|c|c|}
\hline
\texttt{contract} & $\mathcal{O}(L D_{max}^3 d^2)$ \\ \hline
\texttt{multiplyAndOptimize} & $\mathcal{O}(L D_{max}^3 D_A d^2)$ \\ \hline
\texttt{sumAndOptimize} & $\mathcal{O}(L D_{max}^3 d^2)$ \\
\hline
\texttt{multiplyScalar} & $\mathcal{O}(D_{max}^2 d)$ \\
\hline
\texttt{spectralDecomposition} & $\mathcal{O}(K^3)$ \\
\hline
\texttt{checkStop} & $\mathcal{O}(1)$ \\
\hline
\end{tabular}
\label{tab:complexity}
\end{table}

We will conclude this section with an analysis of the complexity of our algorithm. The complexities of the subfunctions of Algorithm \ref{approx_lanczos} are listed in Table \ref{tab:complexity}. For the analysis of \texttt{multiplyAndOptimize}, we have assumed $D_A$ to be of the same order as $D_{max}$. If it were significantly larger, the complexity would change to $\mathcal{O}(L D_{max}^2 D_A^2 d^2)$. Note that this analysis does not extend to the number of sweeps necessary for the optimizations to converge. For the spectral decomposition, we have for simplicity assumed all $T_i$ to be of size $K \times K$. Combining all the different results, we thus find that the overall complexity of Algorithm \ref{approx_lanczos} is in $\mathcal{O}(KL D_{max}^3 D_A d^4)$ with $L = \log(N)$ and since we require $K, D_{max} \in \mathcal{O}(\log N)$, this is equivalent to $\mathcal{O}(poly(L)).$
\section{Numerical Results}
\label{results}
In this section we present numerical results obtained for a challenging problem of relevance in quantum many body physics. Our goal thereby is to study the convergence of the results with increasing $K$ and $D_{max}$. The problem we consider is the approximation of the von Neumann entropy. 
For a quantum state $\rho$~\footnote{$\rho$ is a positive operator with unit trace, representing an ensemble of pure states}, the von Neumann entropy is given by $S = - \trace \rho \log \rho$. In the following, we will focus on the case of states of the form $\rho = e^{-\beta H} \mathbin{/} \trace e^{-\beta H}$, i.e., thermal equilibrium states for a Hamiltonian $H$, at a certain inverse temperature, $\beta$~\footnote{Note that $\beta$ is not related to the $\beta_i$ computed by our algorithm.}. Here, we assume $H$ to be the Ising Hamiltonian with open boundary conditions that is given by 
\[ H = J \sum_{i=1}^{L-1} \sigma_i^x \sigma_{i+1}^x + g \sum_{i=1}^{L} \sigma_i^z + h \sum_{i=1}^{L} \sigma_i^x
\]
where $\sigma^{\{x,y,z\}}$ are the Pauli matrices
\begin{equation*}
\sigma^x = \left[ \begin{matrix}
0 & 1 \\
1 & 0 \\ \end{matrix} \right],
\sigma^y = \left[ \begin{matrix}
0 & -i \\
i & 0 \\ \end{matrix} \right],
\sigma^z = \left[ \begin{matrix}
1 & 0 \\
0 & -1 \\ \end{matrix} \right].
\end{equation*} Note however that here by $\sigma_i^{\{x,y,z\}}$ we actually denote the tensor product $I_1 \dots \otimes I_{i-1} \otimes \sigma_i^{\{x,y,z\}} \otimes I_{i+1} \otimes \dots \otimes I_L$ and analogously for $\sigma_i^{\{x,y,z\}}\sigma_{i+1}^{\{x,y,z\}}$ for simplicity of notation. The Hamiltonian describes a chain of spin particles with nearest neighbour interactions and two magnetic fields acting only on the individual particles. This choice of $H$ has the advantage that it is exactly solvable for $h = 0$, a case commonly known as 'transverse field Ising', and thus opens the possibility to obtain a reference solution for system sizes for which it could otherwise not be obtained \cite{sachdev2007quantum, karevski2006ising, suzuki2012quantum}. 
It is possible to find an MPO-approximation to the thermal equilibrium state, $\rho$, by means of standard MPS-techniques~\cite{verstraete2004matrix,zwolak2004mixed,garcia2006time}. 
It is customary to use a \emph{purification} ansatz for this purpose, where $\rho(\beta/2)$ is approximated by standard imaginary time evolution, and the whole state is then written as $\rho\propto \rho(\beta/2)^* \rho(\beta/2)$. In the context of our algorithm, nevertheless, applying exactly this $\rho$ involves a larger cost and worse numerical condition. Instead, we apply the method as described above to $\rho(\beta/2)$ and absorb the necessary squaring into the function that is to be approximated.
In our case this means that instead of computing $f(\lambda_i) = \lambda_i \log \lambda_i$ for each Ritz-value, we simply compute $f(\lambda_i)^{\prime} = \lambda_i^2 \log \lambda_i^2$. This allows us to apply the algorithm to a possibly much more benign input at the cost of an only slightly more complicated function. 
Due to truncation errors, the operator $\rho(\beta/2)$ may not be exactly Hermitian. This can be easily accounted for by taking its Hermitian part, $\frac{1}{2}[\rho(\beta/2)^*+\rho(\beta/2)]$, which is an MPO with at most twice the original bond dimension.
In our experiments we however did not find this to be necessary.

Apart from the entropy, another interesting function to examine would have been the trace norm given by $|| \rho||_1 = \trace \sqrt{\rho^* \rho}$, i.\ e., the sum of the singular values. But as we only consider positive matrices in this scenario, this sum is equal to the sum of the eigenvalues which we know to be equal to 1 due to the normalization of the thermal state. In fact, we find that for $\rho(\beta)$ $\alpha_1$ as computed by our algorithm is given by 
\[\alpha_1 = \trace U_1^* V_1 = \trace U_1^* \rho U_1 = \frac{1}{\beta_1^2} \trace U^* \rho U = \frac{1}{\beta_1^2} \trace \rho. \] So, in this case the algorithm computes the exact result in one step. We verified this result numerically and found it to hold for all considered cases. It is also easily possible to adapt the above trick to compute the trace norm from $\rho(\frac{\beta}{2})$. We found this to approach to work very well but do not present results here as computing the value directly as shown above is of course preferable. In future work, the trace norm might be used as a distance measure between different thermal states.

As was explained in Section \ref{quadrature}, the $2K$-th derivative of $f$ determines whether $\mathcal{G}f$ poses a lower or upper bound of the true value. In our case and for $K>1$, it is given by
\[\frac{d^{2K} - \lambda_i^2 \ln \lambda_i^2}{d^{2K} \lambda_i} = \frac{4(2K-3)!}{\lambda_i^{2K-2}}\]
with $\lambda_i$ being the $i$-th eigenvalue of $\rho$. Hence, we can expect our algorithm to provide increasingly tight lower bounds for the correct value. We use the violation of this property as a stopping criterion to account for the situation when truncation errors become too large. 
Additionally, we keep the average of the last three or four, depending on $\beta$, approximations and employ the aforementioned 3$\sigma$-rule. In case these stopping criteria are not met, we terminate the algorithm when the difference between successive approximations is below $10^{-10}$.
\begin{figure}[t]
	\begin{subfigure}[t]{0.5\textwidth}
		\centering\captionsetup{width=.8\linewidth}
		\includegraphics[width=1\textwidth]{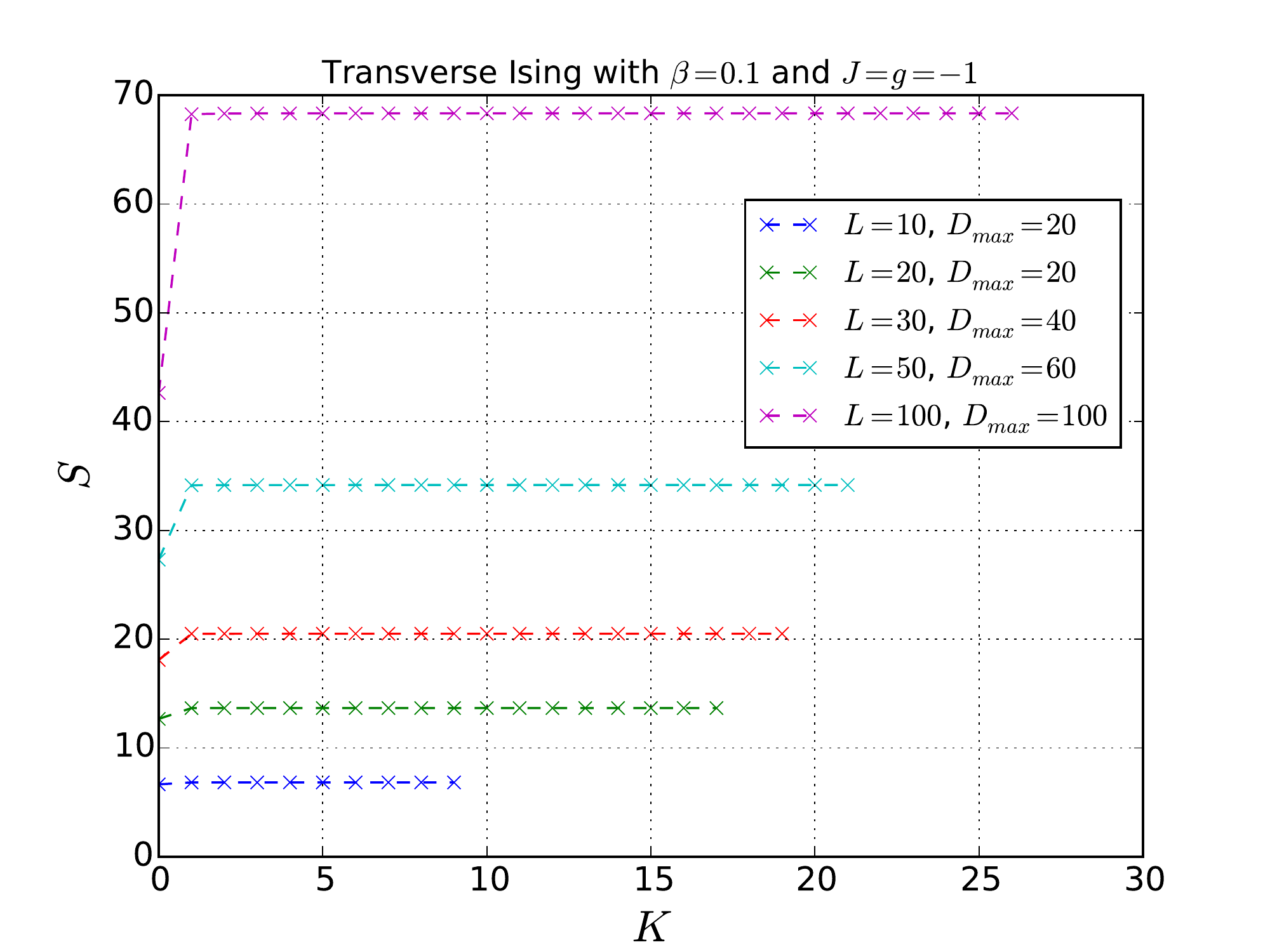}
		\subcaption{Entropy over $K$.}
		\label{result_plots_beta_0_1_conv}
	\end{subfigure}
	\begin{subfigure}[t]{0.5\textwidth}
		\centering\captionsetup{width=.8\linewidth}
		\includegraphics[width=1\textwidth]{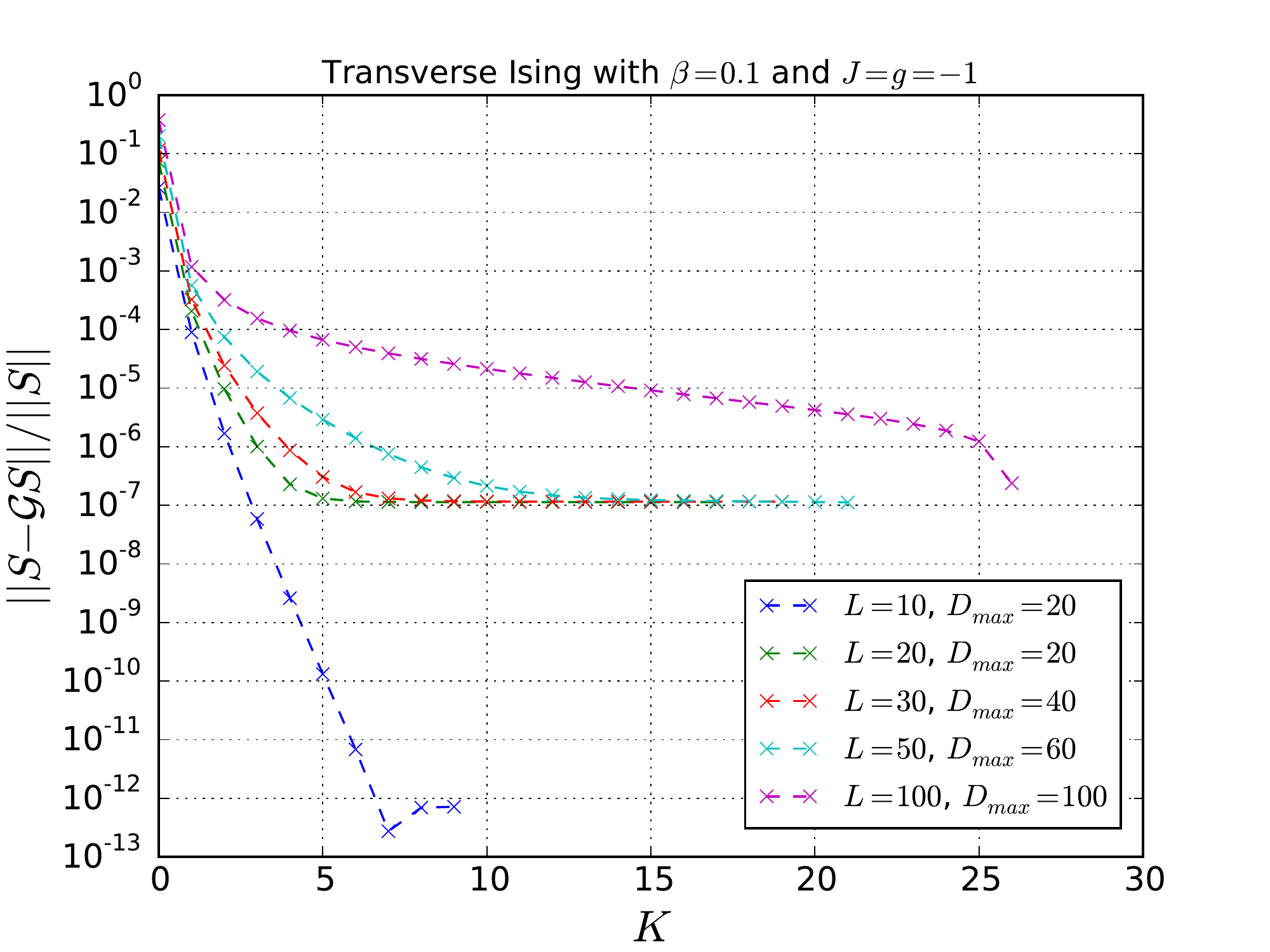}
		\subcaption{Relative error in the entropy over $K$.}
		\label{result_plots_beta_0_1_bond}
	\end{subfigure}
	\caption{Convergence behavior of the algorithm for $L \in \{10,20,30,50,100\}$, $\beta=0.1$ and varying Krylov-dimension $K$. In (a), the convergence of the approximation is depicted. In (b), the convergence of the relative error is shown.}
	\label{beta_0_1_plots}
\end{figure}
\begin{figure}[t]
	\begin{subfigure}[t]{0.5\textwidth}
		\centering\captionsetup{width=.8\linewidth}
		\includegraphics[width=1\textwidth]{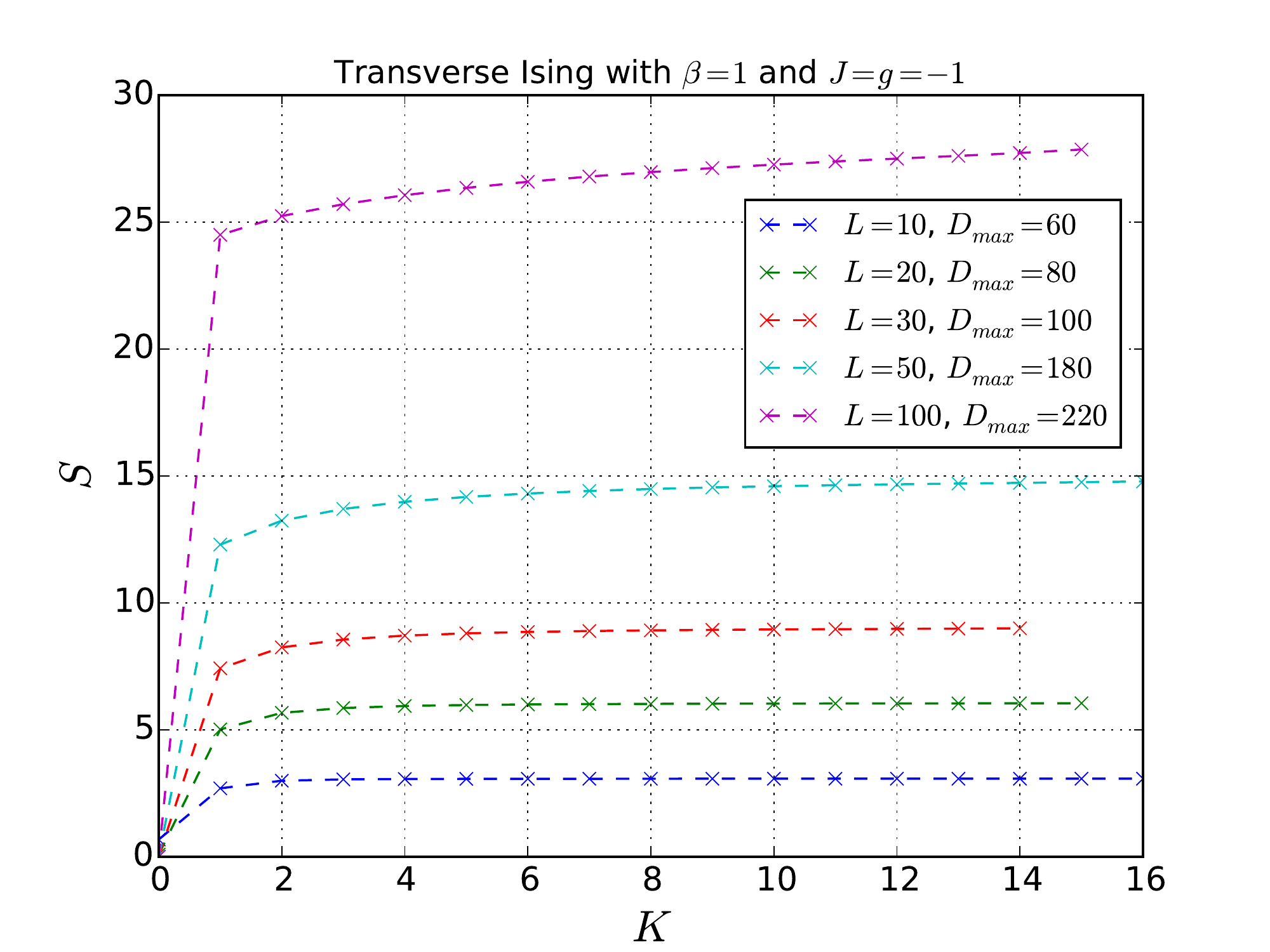}
		\subcaption{Entropy over $K$.}
		\label{result_plots_beta_1_conv}
	\end{subfigure}
	\begin{subfigure}[t]{0.5\textwidth}
		\centering\captionsetup{width=.8\linewidth}
		\includegraphics[width=1\textwidth]{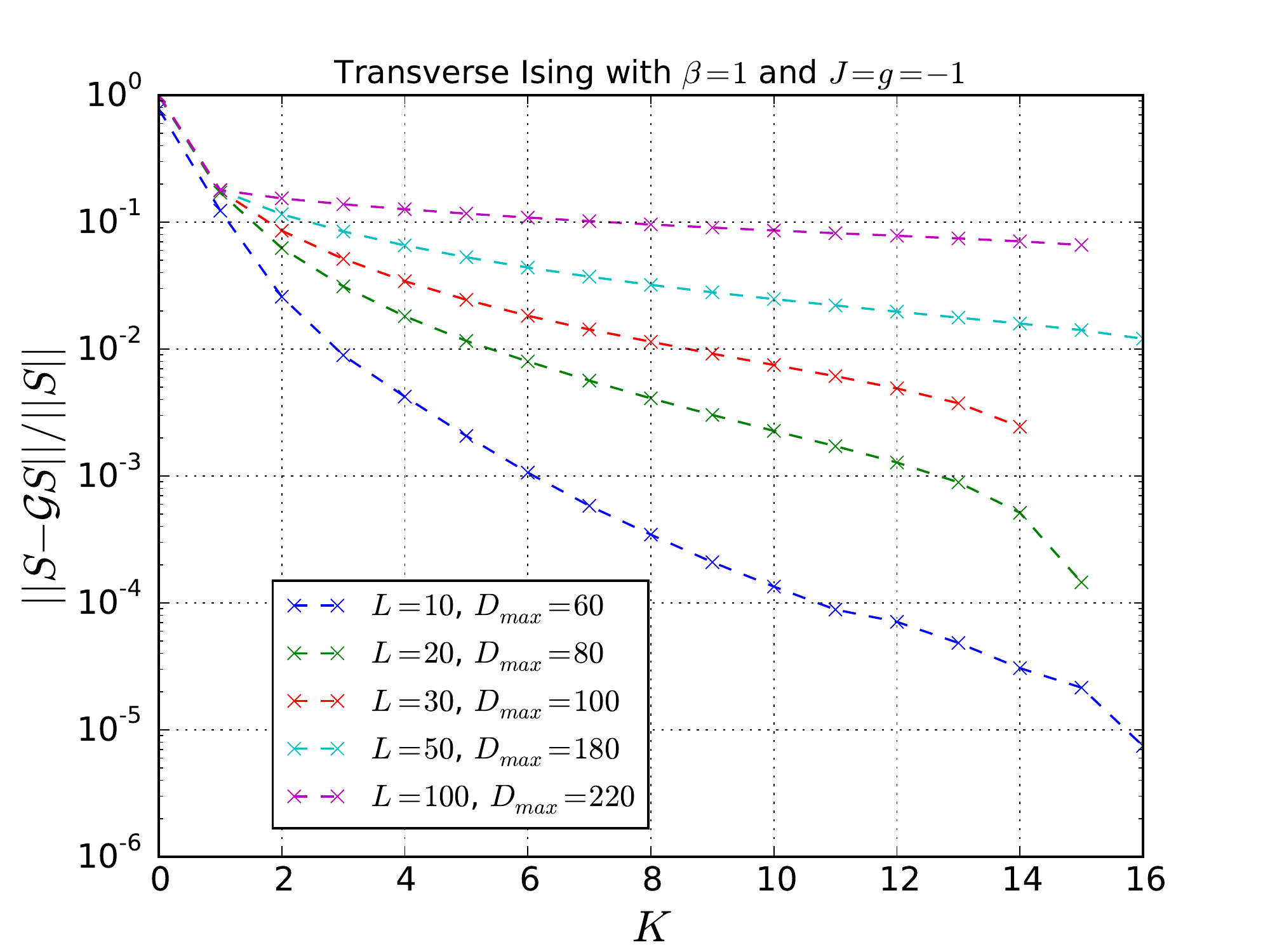}
		\subcaption{Relative error in the entropy over $K$.}
		\label{result_plots_beta_1_bond}
	\end{subfigure}
	\caption{Convergence behavior of the algorithm for $L \in \{10,20,30,50,100\}$, $\beta=1$ and varying Krylov-dimension $K$. In (a), the convergence of the approximation is depicted. In (b), the convergence of the relative error is shown.}
	\label{beta_1_plots}
\end{figure}

\begin{figure}[t]
	\centering\captionsetup{width=.8\linewidth}
	\includegraphics[width=.5\textwidth]{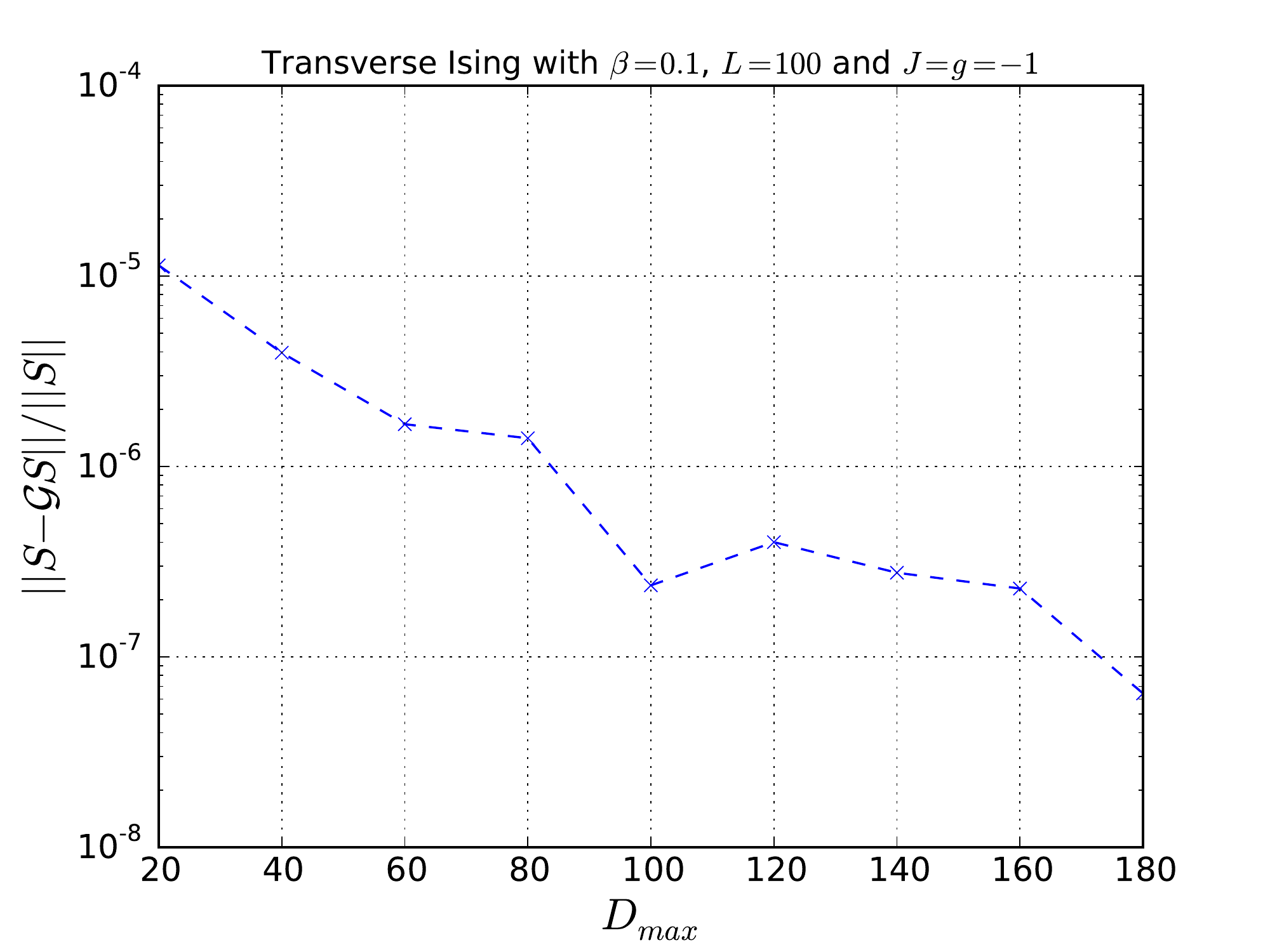}
	\caption{Relative error in the entropy for $L=100$ and $\beta=0.1$ over the maximal bond dimension $D_{max}$. }
	\label{result_plots_100_conv}
\end{figure}

\begin{figure}[t]
	\centering\captionsetup{width=.8\linewidth}
	\includegraphics[width=.5\textwidth]{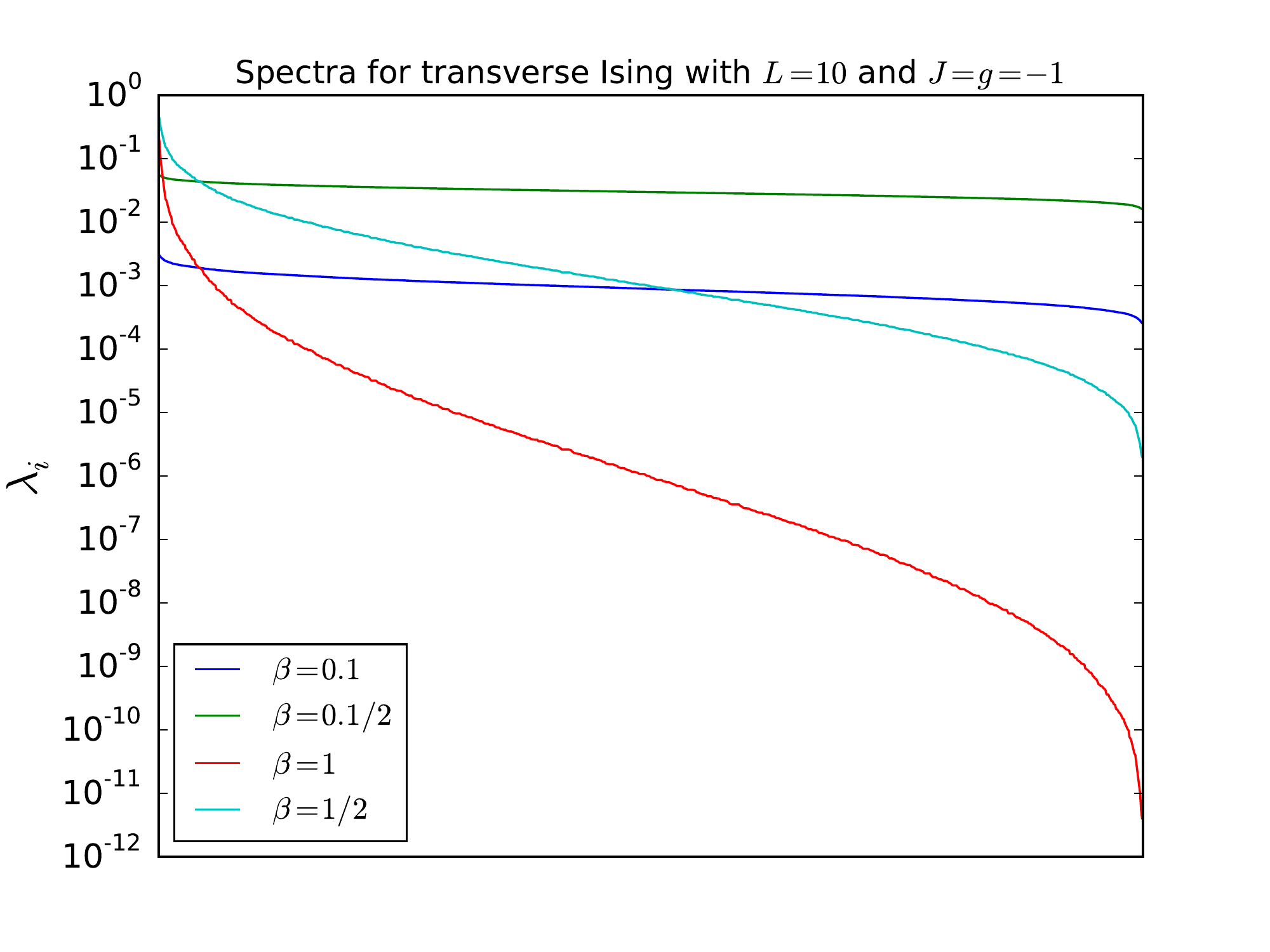}
	\caption{The spectra of $\rho(\beta)$ and $\rho(\frac{\beta}{2})$ for $L=10$. }
	\label{spectra}
\end{figure}
In our experiments we considered systems of size $L \in \{10,20,30,50,100\}$, Hamiltonian parameters $J=g=1$ and inverse temperatures $\beta \in \{0.1, 1.0\}$. The bond dimension used to obtain $\rho(\frac{\beta}{2})$ was set to 20 for all cases. The convergence of the approximation as well as the relative error for $\beta=0.1$ and $\beta=1$ are shown in Figure~\ref{beta_0_1_plots} and Figure~\ref{beta_1_plots}, respectively. Note that in order to compute the relative error, we used numerical diagonalization for $L=10$ and the analytical solution~\cite{sachdev2007quantum} for $L>10$. Figure~\ref{result_plots_100_conv} shows the change of the relative error of the final approximations with growing $D_{max}$ for $L=100$ and $\beta=0.1$.

For the case of $\beta=0.1$ we observe fast convergence to good approximations in $K$ and $D_{max}$, as shown in Figure~\ref{beta_0_1_plots}. The maximal bond dimension required for good convergence only grows mildly with $L$ allowing our method to scale very well with the size of the input. 
The plots in Figure~\ref{result_plots_beta_0_1_bond} show a plateau in the relative error at  $~10^{-7}$. This corresponds to the non-vanishing difference between the exact solution and the numerical MPO used as input. Notice that for $L=10$, where the input is exact, the method is able to achieve a smaller error.
Figure~\ref{result_plots_100_conv} illustrates that while our method achieves a good error even for a small $D_{max}$ of 20, it profits from an increase of the maximal bond dimension. The error decreases by two orders of magnitude from $10^{-5}$ to about $10^{-7}$ when $D_{max}$ is raised from 20 to 180 which still constitutes a strong truncation. It is conceivable that a further increase of the maximal bond dimension would improve the accuracy. We also found that $D_{max}$ limits the number of basis MPOs that can be successfully orthogonalized and therefore effectively controls the maximally reachable $K$. Hence, $D_{max}$ can be regarded as the decisive parameter of our method.

The results for the larger inverse temperature $\beta=1$ paint a slightly different picture. While the overall behavior of our method remains the same and Figure~\ref{result_plots_beta_1_conv} depicts good convergence especially for $L<100$, Figure~\ref{result_plots_beta_1_bond} shows that the relative error achieved is noticeably worse than for the case of $\beta=0.1$. It also seems that larger values of $D_{max}$ are required to achieve reasonable results. This phenomenon naturally becomes more pronounced with larger $L$.

We conjecture that the difference in the performance observed for the two considered values of $\beta$ has two main reasons. Firstly, the bond dimension required for a good approximation of $\rho$ grows with larger $\beta$. This might in turn increase the value of $D_{max}$ required for good accuracy, and, correspondingly, increase the approximation error incurred by $\rho$, so that the computed function will be farther from the analytical solution.
Secondly, the spectral properties of the obtained MPOs for the two considered cases are significantly different. In Figure~\ref{spectra}, we show the spectra of $\rho$ for both values of $\beta$ and $\beta/2$ for $L=10$, respectively. It is clearly visible that $\beta=0.1$ poses a much more benign case. This is underlined by the condition numbers, which are roughly 11.9, $5.68 \cdot 10^{10}$, 3.5 and $2.38 \cdot 10^{5}$ for $\beta=0.1$, $\beta=1$, $\beta=0.05$ and $\beta=0.5$, respectively. They show that $\beta=1$ in fact yields highly ill-conditioned MPOs, functions of which are hard to approximate. These considerations also make it clear that by absorbing the necessary squaring of the eigenvalues into the function, we obtain much more well-conditioned input MPOs of lower bond dimension. Hence, we can conclude that our method, while being influenced by both aforementioned factors, is relatively robust and even works reasonably well for very difficult cases.

We do not provide a comparison to other methods at this point, because of the simple reason that to the best of the authors knowledge there is no other algorithm that can solve the considered kind of problem for $L \gg 20$. For $L \leq 20$ we however expect the existing highly optimized methods to outperform our method in terms of runtime.

\section{Discussion}
\label{discussion}
In this work, we have introduced a method to approximate functionals of the form $\trace f(A)$ for matrices of dimension much larger than $2^{20}$. We started by giving an overview over the mathematical and algorithmic ideas behind the method. Following this, a detailed description of the algorithm together with an analysis of its complexity was provided. We then presented numerical results for a challenging problem in quantum many body physics. These results indicate that our method is able to produce good approximations for a number of Krylov steps and a maximal bond dimension logarithmic in the size of the matrix as long as the matrix exhibits some structure that can be expressed well in the MPO/MPS-formalism and is moderately well-conditioned. It was also shown that the maximal allowed bond dimension is the decisive parameter of the algorithm.

There are several ways to build upon this work. Firstly, an investigation of preconditioning methods suitable for our method could be fruitful. Secondly, a more thorough analysis of the effect of the approximation error introduced by the tensor network formalism on the approximation error of the Gauss quadrature would be an interesting addition. Thirdly, the connection of the approximability of a matrix by an MPO to the convergence behavior of our method could provide deeper understanding. Fourthly, it could be investigated which of the many improvements over the normal Gauss quadrature, as for instance \cite{reichel2015generalized}, can be incorporated into our algorithm to make better use of the expensive information obtained in the Krylov iteration. Finally, the method naturally could be applied to solve practical problems of interest. 

While our method was tested for a quantum mechanical problem, it is of course general in nature and can be applied to any case where the matrix in question can be formulated as an MPO or well approximated by one. Especially for matrices of dimension larger than $2^{10}$ that however can still be explicitly stored, it might be interesting to consider computing the desired function for the MPO-representation.

\section*{Acknowledgements}
This work was partly funded by the \emph{Elite Network of Bavaria}(ENB) via the doctoral programme \emph{Exploring Quantum Matter}.



\nocite{*}
\bibliographystyle{alpha}
\bibliography{ref}

\end{document}